\documentclass[5p]{elsarticle}
\usepackage{amssymb, amsmath, psfrag, ae, graphicx, longtable, array}
\usepackage{algorithmic, algorithm, epstopdf, multirow, mathtools}
\usepackage{color,float,caption,subcaption}
\usepackage{mathcomp,textcomp}
\definecolor{mygreen}{RGB}{0, 150, 0}

\usepackage{units,eurosym,booktabs}
\def\ve#1{\mathchoice{\mbox{\boldmath$\displaystyle#1$}}
{\mbox{\boldmath$\textstyle#1$}}
{\mbox{\boldmath$\scriptstyle#1$}}
{\mbox{\boldmath$\scriptscriptstyle#1$}}} 
\usepackage{lscape}
\usepackage{soul}
\usepackage{color}

\newcommand{\dd}{\text{d}}
\usepackage{tikz}
\newcommand\Square[1]{+(-2.5pt,-2.5pt) rectangle +(2.5pt,2.5pt)}
\newcommand\Rrectangle[1]{+(-7.5pt,-3pt) rectangle +(7.5pt,3pt)}
\newcommand\LongRrectangle[1]{+(-4pt,-0.1pt) rectangle +(4pt,0.1pt)}
\usepackage{pgfplots}
\usetikzlibrary{plotmarks}
\usepackage{xspace} 
\usepackage[
textsize=scriptsize]{todonotes}

\newcommand{\argmin}[1]{\underset{#1}{\operatorname {arg} \! \operatorname {min}}\;}
\newcommand{\argmax}[1]{\underset{#1}{\operatorname {arg} \! \operatorname {max}}\;}

\makeatletter
\newcommand*{\rom}[1]{\expandafter\@slowromancap\romannumeral #1@}
\makeatother
    
\journal{Journal of Process Control}

\begin{document}
\begin{frontmatter}

\title{Optimal Experiment Design in Nonlinear Parameter Estimation with Exact Confidence Regions}

\author[TUDO]{Anwesh Reddy Gottu Mukkula}
\ead{anweshreddy.gottumukkula@tu-dortmund.de}
\author[TUDO,STUBA]{Radoslav Paulen\corref{cor1}}
\ead{radoslav.paulen@stuba.sk}
\cortext[cor1]{Tel.: +421 (0)2 5932 5730, Fax: +421 (0)2 5932 5340 (R. Paulen)}
\address[TUDO]{Process Dynamics and Operations Group, Department of Chemical and
Biochemical Engineering, Technische Universit\"at Dortmund, Emil-Figge-Strasse
70, 44227 Dortmund, Germany}
\address[STUBA]{Faculty of Chemical and Food Technology, Slovak University of 
Technology in Bratislava, Radlinsk\'eho 9, 812 37 Bratislava, Slovakia}

\begin{abstract}
A model-based optimal experiment design (OED) of nonlinear systems is studied. OED represents a methodology for optimizing the geometry of the parametric joint-confidence regions (CRs), which are obtained in an a posteriori analysis of the least-squares parameter estimates. The optimal design is achieved by using the available (experimental) degrees of freedom such that more informative measurements are obtained. Unlike the commonly used approaches, which base the OED procedure upon the linearized CRs, we explore a path where we explicitly consider the exact CRs in the OED framework. We use a methodology for a finite parametrization of the exact CRs within the OED problem and we introduce a novel appro\-xi\-mation technique of the exact CRs using inner- and outer-appro\-xi\-mating ellipsoids as a computationally less demanding alternative. The employed techniques give the OED problem as a finite-dimensional mathematical program of bilevel nature. We use two small-scale illustrative case studies to study various OED criteria and compare the resulting optimal designs with the commonly used linearization-based approach. We also assess the performance of two simple heuristic numerical schemes for bilevel optimization within the studied problems.
\end{abstract}

\end{frontmatter}

\section{Introduction}\label{sec:intro}
At present, advanced industrial engineering and management strive for resource- and energy-efficient design and operation of systems, plants, and processes. Here a use of the model-based techniques is a leading paradigm. The employed models, whether mechanistic or data-based, include a finite number of parameters, whose values are related to the particular natural and system-wide phenomena and are thus commonly only known to belong to some interval or unknown completely. Therefore, as-precise-as-possible determination of the unknown (uncertain) model parameters is crucial for the deployment of the effective model-based solutions.

In the world, where the sensing technology becomes virtually present everywhere, the simplest way of obtaining the parameter values is to conduct a series of observations, experiments, during which the measurements of some quantities (output variables) are gathered. Subsequently, an estimation procedure is employed to find parameter values such that the model outputs match the observed data. As the measurements from a real plant are corrupted with some (random) noise, the uncertain model parameters cannot be identified precisely. This a posteriori uncertainty can be represented by a joint-confidence region (CR) that encompasses all the likely parameter estimates, given the probability distribution of the measurement noise.

The parametric uncertainty can be reduced by performing an experiment that appropriately sets the plant into an operating region, where more informative measurements are obtained. A way how to identify such experiment consists in performing an optimal experiment design~(OED)~\citep{hja05, pro08}. OED identifies the best experiment in terms of initial conditions, control inputs, sampling times and/or locations of measurement devices. The model-based OED problem is usually formulated as a mathematical program, where a certain measure of the CR, such as volume, is minimized.

Some well-established methods exist for the OED for linear systems~\citep{fra08}, wherein CR is an ellipsoid~\citep{SeberWild200309}, so the task of its optimal shaping is greatly simplified. For nonlinear systems, the most common, yet appro\-xi\-mate, approach is to resort to a system linearization and a use of the linear techniques. Beale~\citep{beale1960} presented methodology for assessment of system nonlinearity in this respect. Other approaches, covered in~\cite{pro13} {\color{black}for convex problems}, are of more or less appro\-xi\-mate character and mostly rely on the {\color{black}convexity} properties of the OED problem. Bates and Watts~\citep{bates1988} presented a framework based on local nonlinear curvature properties, which is also an appro\-xi\-mate technique. The use of a finite parametrization of the exact CRs~\citep{kur97, rooney2001, str14} is a relatively recent subject of study.

In this contribution, we study the framework for OED of nonlinear systems that is based on explicit consideration of the exact CRs. We also present a computationally simpler variant that is based on simultaneous inner- and outer-appro\-xi\-mations of the CR by ellipsoids. The preliminary findings of this work were presented in the conference contribution~\citep{mupa_ifac17}. We study various common design criteria and compare the performance of the presented techniques with the linearized OED.

We organized the paper as follows. The concepts of linear and nonlinear parameter estimation and construction of CRs are discussed first. Next, the formulation is presented of the experiment design criteria that directly use exact CRs (exact designs). We also present {\color{black}two simple heuristic} numerical approaches taken from literature to solve the arising bilevel optimization problems. Finally results of two {\color{black}illustrative} case studies are presented and discussed.

\section{Parameter Estimation Problem}
\subsection{Mathematical Model}
In this paper, a mathematical model of a system is represented by
\begin{equation}
  \label{eq:model_static}
    \quad \hat{\ve y}(\ve p, \tau) = \ve F(\ve p, \ve u_\tau),
\end{equation}
with $\hat{\ve{y}}$ as $n_y$ measured variables, $\ve u_\tau$ as $n_u$ degrees of freedom and $\hat{\ve{p}}$ as $n_p$ uncertain parameters. Here $\tau$ represents an ordinal number of the data point taken in one or more experiments. Function $\ve F:\mathbb{R}^{n_p}\times\mathbb{R}^{n_u}\rightarrow\mathbb{R}^{n_y}$ {\color{black}is} {\color{black}twice continuously differentiable} mapping. Throughout the paper, we resort to a representation~\eqref{eq:model_static}, which considers the system model as static and explicit. However, the presented findings can {\color{black}straightforwardly }be extended to dynamic and implicit models.

We will assume throughout the paper that the model is not over-parametrized and that the parameters are identifiable. We consider that, upon the realization of an experiment or several experiments, $N$ instances are gathered of $n_y$-dimensional vector of plant measurements $\ve{y}_m$ and are subsequently used for the estimation of unknown parameters. Throughout the paper, we will assume the Gaussian noise to be corrupting the measurements. In the following subsections, existing frameworks are presented for the identification of the unknown parameters and the corresponding exact CRs for both the linear and the nonlinear parameter estimation problems.

\subsection{Linear Parameter Estimation}\label{subsec:Linear_Parameter_Estimation}
Parameter estimation with a linear model involved is a well-studied topic in the literature~\citep{SeberWild200309}. Assume a mathematical model with a mapping function $\ve F_{l}$ of the form
\begin{align} \label{eq:linear}
  \hat{\ve{y}}(\ve{p},\tau) = \ve{F}_{l}(\ve{p},\ve{u}_{\tau})
  = \mathbf{Q}(\ve u_\tau) \ve p,
\end{align}
where $\mathbf{Q}(\ve u_\tau)$ is a so-called regressors matrix with appropriate dimensions.

Under the assumption of uncorrelated and normally distributed measurement noise with \emph{known standard deviation} vector $\ve\sigma$, the maximum-likelihood estimate is found via the least-squares estimation as
\begin{align}
\hat{\ve p} &\in \argmin{\ve{p}} J_w(\ve p), \label{eq:w_with_variance}
\end{align}
{\color{black}with}
\begin{align}
J_w(\ve p)&:= \sum_{i=1}^{n_y} \sum_{\tau=\tau_1}^{\tau_N} \sigma_i^{-2} (y_{i,m}(\tau) - \hat{{y}}_{i}(\ve{p},\tau))^2.
\end{align}
The CR of parameter estimates is then
given by an ellipsoid~\citep{SeberWild200309}
\begin{equation}\label{eq:LPE}
(\ve{p} - \hat{\ve{p}})^T \mathbb{FIM}(\mathcal U)
(\ve{p} - \hat{\ve{p}}) \leq \chi^2_{\alpha,n_p},
\end{equation}
where $\mathbb{FIM}$ is the so-called Fisher information matrix,
\begin{equation}
\mathbb{FIM}(\mathcal U) = \sum_{\tau=\tau_1}^{\tau_N} \mathbf Q(\ve u_\tau)^T \text{diag}^{-1}(\sigma_1^2,\hdots,\sigma_{n_y}^2) \mathbf Q(\ve u_\tau),
\end{equation}
$\mathcal U:=(\ve u_{\tau_1}^T,\ve u_{\tau_2}^T,\dots,\ve u_{\tau_N}^T)^T$, and $\chi^2_{\alpha, n_p}$ represents the upper $\alpha$ quantile of the chi-squared statistical distribution with $n_p$ degrees of freedom.
%
%

If the \emph{variance of the measurement noise is unknown}, the covariance matrix $\text{diag}(\sigma_1^2,\hdots,\sigma_{n_y}^2)$ is normally appro\-xi\-mated by $s^2\mathbb I$ with
\begin{align} \label{eq:variance}
  s^2 := \frac{J(\hat{\ve{p}})}{N-n_p},
\end{align}
where the expected (most-likely) value of parameters~$\hat{\ve p}$ is identified by solving
\begin{equation}\label{eq:residual}
\hat{\ve p} \in \argmin{\ve{p}} J(\ve p) = \sum_{i=1}^{n_y} \sum_{\tau=\tau_1}^{\tau_N} (y_{i,m}(\tau) - \hat{{y}}_{i}({\ve p},\tau))^2.
\end{equation}
The joint-confidence ellipsoid is then given by
\begin{align}
  (\ve{p} - \hat{\ve{p}})^T \left(\sum_{\tau=\tau_1}^{\tau_N}\mathbf Q(\ve u_{\tau})^T \left(s^2\mathbb I\right)^{-1} \mathbf Q(\ve u_{\tau})\right)
(\ve{p} - \hat{\ve{p}}) \notag \\ \leq n_p \mathcal F_{n_p,N-n_p,\alpha},
\label{eq:LPE_unknown_variance}
\end{align}
where $\mathbb{FIM}$ is replaced by its corresponding appro\-xi\-mation and $\mathcal F$ represents the upper $\alpha$ quantile of the Fisher distribution with $n_p$ and $N-n_p$ degrees of freedom in the numerator and in the denominator, respectively.


\subsection{Nonlinear Parameter Estimation}
Given a static nonlinear mathematical model
\begin{align} \label{eq:nonlinear}
  \hat{\ve{y}}(\ve{p},\tau) = \ve{F}_{nl}(\ve{p},\ve{u}_\tau),
\end{align}
one can identify the (exact) CR dependent upon the availability of information about the variance of the measurement noise. If the variance is known, the exact CR is given by~\citep{SeberWild200309}
\begin{align} \label{eq:NLPE_with_variance}
  J_w(\ve{p}) - J_w(\hat{\ve{p}}) \leq \chi^2_{\alpha,n_p},
\end{align}
while, when the variance of the measurement noise is unknown, the exact CR is given by~\citep{SeberWild200309}
\begin{align} \label{eq:NLPE_without_variance}
  J(\ve{p}) - J(\hat{\ve{p}}) \leq n_p s^2 \mathcal F_{n_p,N-n_p,\alpha}.
\end{align}
At this point we can define the sets $P_w:=\{\ve p \,|\, \text{Eq.}~\eqref{eq:NLPE_with_variance}\}$ and $P:=\{\ve p \,|\, \text{Eq.}~\eqref{eq:NLPE_without_variance}\}$. Unlike in the linear parameter estimation, the CR does not generally have a shape of an ellipsoid due to nonlinearity.

\section{Optimal Experiment Design}\label{sec:OED}
We present a methodology for OED for both linear and nonlinear parameter estimation problems.
We will assume that the CR is connected. For disjoint exact CRs, which result from non-identifiability issues, it is normally suggested to perform a re-parameteri\-za\-tion of the model~\citep{bates1988}. We will also assume that an estimate $\hat{\ve p}$ is available. The final assumptions which are inherent to the standard experiment design techniques is that there exists no structural plant-model mismatch and that the expected realization of the measurement noise is 0. In turn this results in $\ve y_m(\tau)=\hat{\ve y}(\hat{\ve p}, \tau),\,\forall \tau$.

Several design criteria are proposed in the literature~\citep{fra08} such as A, D, E, Modified E, V, Q, M and so on. Each of these designs aims to tune a particular property of the confidence region. We will focus our study on the most used criteria, i.e., A, D, and E, but other design criteria might be considered as well using the ideas presented herein.

\subsection{A-optimal design}
\begin{figure}[tb!]
  \centering
    \psfrag{P1}[cc][cc][1]{$p_1$}
    \psfrag{P2}[cc][cc][1]{\ $p_2$}
    \psfrag{True E design}[cc][cc][1]{\quad E design}
    \psfrag{True A design}[cc][cc][1]{\quad A design}
    \psfrag{True D design}[cc][cc][1]{\quad D design}
    \centering\includegraphics[width=0.95\linewidth]{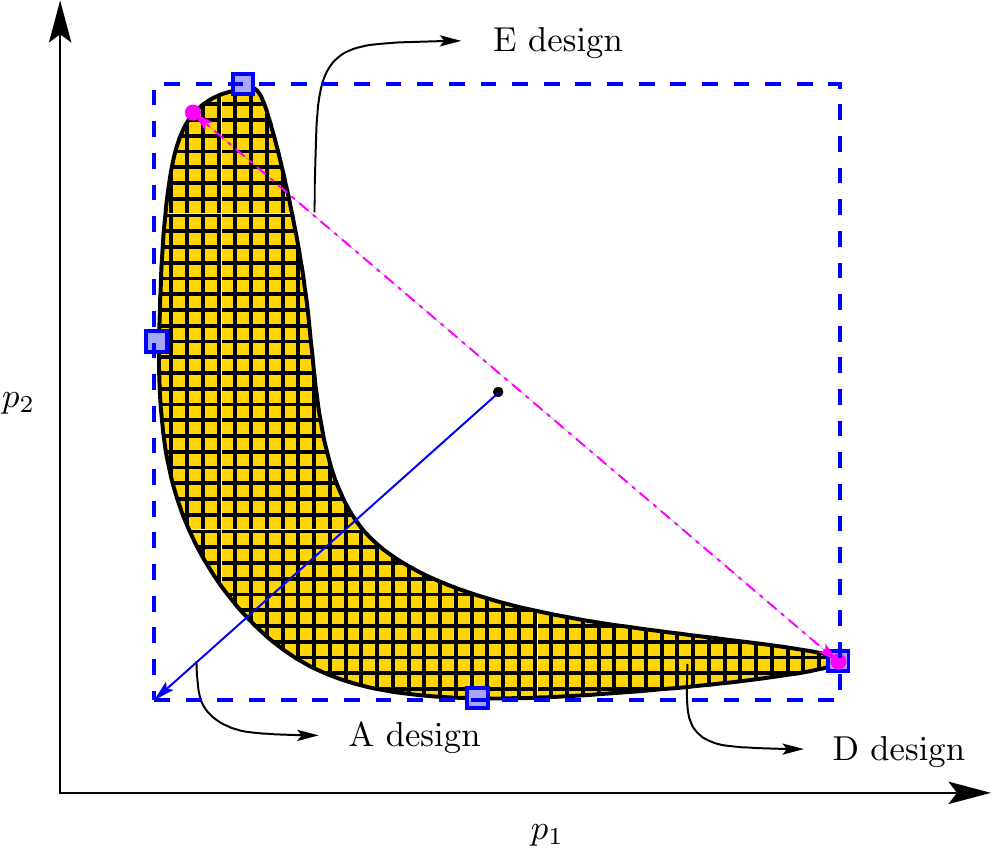}
    \caption[]{Illustration of the design criteria in two-dimensional parametric space. The plot shows a generic exact CR (shaded), the over-approximating orthotope (dashed) of the exact CR identified using the anchor points $\ve \pi$ (\tikz\draw[line width=0.4 mm,blue, fill={rgb,255:red,180; green,170; blue,255}] (0,0) \Square{10pt};). \tikz\draw[line width=0.3 mm, white,fill={rgb,255:red,255; green,100; blue,255}] (0,0) circle (.75ex); mark the points that give the maximal Euclidean distance of two points in the CR.}
\label{fig:designs}
\end{figure}

The idea behind the A design criterion is to minimize the perimeter of the box that encloses the exact CR~\citep{fra08}, i.e., to minimize the sum of projections of the CR on the parameter axes. This idea is sketched in Fig.~\ref{fig:designs}, where the shaded set represents a CR. The enclosing box is given by the four anchor points, marked as 
squares (\tikz\draw[line width=0.4 mm,blue, fill={rgb,255:red,180; green,170; blue,255}] (0,0) \Square{10pt};).

For a general CR, one can identify $2n_p$ anchor points
\begin{align}\label{eq:pi}
	\ve\pi\!:=\!
	\left\{\!\!
	\left(\!\!\begin{array}{c}
		p_1^L\\ p_2^{1,L}\\ \vdots\\ p_{n_p}^{1,L}
	\end{array}\!\!\right)\!\!,\!
	\left(\!\!\begin{array}{c}
		p_1^U\\ p_2^{1,U}\\ \vdots \\ p_{n_p}^{1,U}
	\end{array}\!\!\right)\!\!,
	\left(\!\!\begin{array}{c}
	p_1^{2,L}\\ p_2^L \\ \vdots\\ p_{n_p}^{2,L}
	\end{array}\!\!\right)\!,
	\cdots\!\:,\!
	\left(\!\!\begin{array}{c}
		p_1^{n_p,U}\\ p_2^{n_p,U}\\ \vdots\\ p_{n_p}^{U}
	\end{array}\!\!\right)\!\!\right\},
\end{align}
where each point gives a lower or an upper limit of the value of a particular parameter in the exact CR. The anchor points can be identified by solving the following optimization problem
\begin{subequations}\label{eq:oa_ortho_with}
  \begin{align}
  \phi_A(\mathcal U):=\max_{\ve\pi} & \sum_{j=1}^{n_p}p_j^U-p_j^L\\
  \text{s.t. } & \forall j\in\{1,\dots,2n_p\}, \ \forall \tau\in\{\tau_1,\dots,\tau_N\}:\\
  &\hat{\ve y}(\ve \pi_j, \tau) = \ve F_{nl}(\ve \pi_j, \ve u_{\tau}),\\
  &\ve y_m(\tau)=\hat{\ve y}(\hat{\ve p}, \tau),\\
  \label{eq:oed_nonlinear_ortho_CR1}
  &J_w(\ve\pi_j) - J_w(\hat{\ve{p}}) \leq \chi^2_{\alpha,n_p}.
  \end{align}
\end{subequations}
Note that the problem of identifying the orthotopic enclosure of the CR is formulated for the case when the measurement-noise variance is known using the expression for the CR~\eqref{eq:NLPE_with_variance}. This would simply be exchanged with the expression~\eqref{eq:NLPE_without_variance} in case when the variance is unknown. 

Note that the problem~\eqref{eq:oa_ortho_with} is separable and highly structured. On the other hand it is non-convex in general and the number of its optimization variables ($2n_p^2$) grows quadratically with the number of parameters. This means that identification of an orthotope might get challenging for the state-of-the-art solvers and for the high-dimensional problems.

The A-optimal experiment design (A design) can be identified by solving
\begin{equation}\label{eq:oed_nonlin_A}
  \min_{\ve u_\tau \in [\ve u^L, \ve u^U], \,\forall \tau\in\{\tau_1,\dots,\tau_N\}} \phi_A(\mathcal U),
\end{equation}
which is a special case of a bilevel program. The bounds $\ve u^L$ and $\ve u^U$ represent the lower and upper limits of the experimental degrees of freedom.

In the linear case, the A design is identified by
\begin{subequations}\label{eq:oed_lin_A}
\begin{align}
  &\min_{\stackrel{\ve u_\tau \in [\ve u^L, \ve u^U]}{\forall \tau\in\{\tau_1,\dots,\tau_N\}}} \tilde\phi_\text{A}(\mathbb{FIM}):=\min_{\stackrel{\ve u_\tau \in [\ve u^L, \ve u^U]}{\forall \tau\in\{\tau_1,\dots,\tau_N\}}} \text{trace}(\mathbb{FIM}^{-1})\\
  &\qquad\text{s.t. } \hat{\ve{y}}(\ve{p},\tau) = \ve{F}_{l}(\ve{p},\ve{u}_{\tau}),\ \forall \tau \in\{\tau_1,\dots,\tau_N\}.
\end{align}
\end{subequations}
In case of unknown variance of measurement noise, the appro\-xi\-mate $\mathbb{FIM}$ from~\eqref{eq:LPE_unknown_variance} will be used. We denote this lineari\-zation-based approach as \emph{classical} in this study.

\subsection{D-optimal design}\label{subsec:trueDOED}
The D design aims to find the experimental conditions such that the exact CR would have minimum volume. The D-optimal design problem then reads as
\begin{equation}\label{eq:}
 \min_{\stackrel{\ve u_\tau \in [\ve u^L, \ve u^U]}{\forall \tau\in\{\tau_1,\dots,\tau_N\}}} \phi_D (\mathcal U) := \min_{\stackrel{\ve u_\tau \in [\ve u^L, \ve u^U]}{\forall \tau\in\{\tau_1,\dots,\tau_N\}}} \idotsint\limits_{P_w(\mathcal U)} \dd p_1\dots\dd p_{n_p}.
\end{equation}
As there is no finite-dimensional parameterization of the set $P_w(\mathcal U)$ available in general, it is very hard in general to evaluate the volume integral. We propose to use a gridding-based approach, where the grid is evaluated inside the aforementioned orthotopic enclosure of the CR. The proposed optimization problem may be written as
\begin{subequations}\label{eq:oed_nonlin_D}
  \begin{align}
  &\min_{\stackrel{\ve u_\tau \in [\ve u^L, \ve u^U]}{\forall \tau\in\{\tau_1,\dots,\tau_N\}}} \hat\phi_D(\mathcal U):=\min_{\stackrel{\ve u_\tau \in [\ve u^L, \ve u^U]}{\forall \tau\in\{\tau_1,\dots,\tau_N\}}} \sum_{\forall i\in\mathcal I_\Pi} \delta_i\\
  &\text{s.t. } \delta_i = \begin{cases}
                           1, & \text{if } \Pi_i\in P_w\\
                           0, & \text{otherwise}
                          \end{cases}\label{eq:ind_vol}\\
   & \qquad \Pi = \{p_1^L, p_1^L+\epsilon,p_1^L+2\epsilon,\dots,p_1^U\}\times\notag\\
   & \qquad \quad\{p_2^L, p_2^L+\epsilon,p_2^L+2\epsilon,\dots,p_2^U\} \times\dots\notag\\
   & \qquad \quad\times\{p_{n_p}^L, p_{n_p}^L+\epsilon,p_{n_p}^L+2\epsilon,\dots,p_{n_p}^U\},\\
  & \qquad \max_{\ve\pi} \sum_{j=1}^{n_p}p_j^U-p_j^L,\\
  & \qquad \quad\text{s.t. } \forall j\in\{1,\dots,2n_p\}, \ \forall \tau\in\{\tau_1,\dots,\tau_N\}:\\
  &\qquad \quad\quad\hat{\ve y}(\ve \pi_j, \tau) = \ve F_{nl}(\ve \pi_j, \ve u_{\tau}),\\
  &\qquad \quad\quad\ve y_m(\tau)=\hat{\ve y}(\hat{\ve p}, \tau),\\  
  &\qquad \quad\quad J_w(\ve\pi_j) - J_w(\hat{\ve{p}}) \leq \chi^2_{\alpha,n_p},
  \end{align}
\end{subequations}
where $\epsilon>0$ is the tuning parameter that determines the accuracy of the appro\-xi\-mation and $\mathcal I_\Pi$ is the index set of $\Pi$. This approach for appro\-xi\-mating the volume of the CR is illustrated in Fig.~\ref{fig:designs} as a grid in the shaded set.

In principle, the identification of the orthotopic enclosure of the CR can be removed and the problem~\eqref{eq:oed_nonlin_D} can be modified to a single-level mathematical program. Nonetheless the optimization problem~\eqref{eq:oed_nonlin_D} is non-smooth due to the presence of indicator function~\eqref{eq:ind_vol} and thus it can get very challenging and computationally highly expensive, especially in higher dimensions. An alternative approach can be exploited by appro\-xi\-mation of the volume using the semi-algebraic sets~\citep{DABBENE2017110}. Approaches to inner approximation of the CR based on an orthotope and on successive SDP approximations are presented by Streif et al.~\citep{STREIF2013321}. We propose a simpler appro\-xi\-mation here, which uses an idea similar to the L\"owner-John ellipsoids~\citep{Ball1992}. We construct the inner- and outer-appro\-xi\-mating ellipsoids, which are the scaled counterparts of the linearized CRs. The proposed appro\-xi\-mate D design is found by
\begin{subequations}\label{eq:oed_nonlin_D_app}
  \begin{align}
  &\min_{\stackrel{\ve u_\tau \in [\ve u^L, \ve u^U]}{\forall \tau\in\{\tau_1,\dots,\tau_N\}}} \frac{\text{det}(k_{out}\mathbb{FIM}^{-1})}{k_{out}} + \frac{\text{det}(k_{in}\mathbb{FIM}^{-1})}{k_{in}}\\ \label{subeq:inner_outer_approx_start}
  &\text{s.t. } \max_{\ve p_{out}, \ve p_{in}, k_{out}, k_{in}} k_{out} - k_{in}\\ 
  & \qquad \quad\text{s.t. } \forall \tau\in\{\tau_1,\dots,\tau_N\}:\\
  &\qquad \quad\quad\hat{\ve y}(\ve p_{out}, \tau) = \ve F_{nl}(\ve p_{out}, \ve u_{\tau}),\\
  &\qquad \quad\quad\hat{\ve y}(\ve p_{in}, \tau) = \ve F_{nl}(\ve p_{in}, \ve u_{\tau}),\\
  &\qquad \quad\quad\ve y_m(\tau)=\hat{\ve y}(\hat{\ve p}, \tau),\\
  &\qquad \quad\quad J_w(\ve p_{out}) - J_w(\hat{\ve{p}}) \leq \chi^2_{\alpha,n_p},\\
  &\qquad \quad\quad J_w(\ve p_{in}) - J_w(\hat{\ve{p}}) \leq \chi^2_{\alpha,n_p},\\
  &\qquad \quad\quad (\ve p_{out} - \hat{\ve{p}})^T \mathbb{FIM}(\ve p_{out} - \hat{\ve{p}}) \leq k_{out},\\
  &\qquad \quad\quad (\ve p_{in} - \hat{\ve{p}})^T \mathbb{FIM}(\ve p_{in} - \hat{\ve{p}}) \geq k_{in},\label{subeq:inner_outer_approx_end}  
  \end{align}
\end{subequations}
where $\ve p_{out}$ and $\ve p_{in}$ are intersection points between outer-/inner-appro\-xi\-mating ellipsoids and the exact CR. The scaling factors $k_{out}$ and $k_{in}$ express the magnitude of deviation of the outer- and, respectively, inner-appro\-xi\-mating ellipsoid from the linearized CR. The weighting in the cost function is then introduced to penalize the contribution of the most deviating ellipsoid. This prevents the design procedure from concentrating on shaping the ellipsoid that is potentially a very loose appro\-xi\-mation of the CR and in practice avoids numerical and irregularity problems. Hence that the proposed problem also scales well, i.e., linearly, w.r.t. the number of the parameters as the lower-level problem optimizes $2n_p+2$ variables. The proposed approximate D-optimal design is therefore computationally a less expensive problem when compared to the exact D-optimal design using the gridding-based approach. We will denote the proposed approximate D-optimal design approach as the \emph{ellipsoidal D design}. The idea behind this method, slightly modified, could also be used for an appro\-xi\-mate A design but we do not explore this path in the present study explicitly.

Note also that, if the CR can be expressed exactly as~\eqref{eq:LPE}, the proposed optimization problem boils down to the classical D design where one solves
\begin{subequations}\label{eq:oed_lin_D}
\begin{align}
  &\min_{\stackrel{\ve u_\tau \in [\ve u^L, \ve u^U]}{\forall \tau\in\{\tau_1,\dots,\tau_N\}}} \tilde\phi_\text{D}(\mathbb{FIM}):=\min_{\stackrel{\ve u_\tau \in [\ve u^L, \ve u^U]}{\forall \tau\in\{\tau_1,\dots,\tau_N\}}} \text{det}(\mathbb{FIM}^{-1})\\
  &\qquad\text{s.t. } \hat{\ve{y}}(\ve{p},\tau) = \ve{F}_{l}(\ve{p},\ve{u}_{\tau}),\ \forall \tau \in\{\tau_1,\dots,\tau_N\}.
\end{align}
\end{subequations}

\subsection{E-optimal design}\label{subsec:trueEOED}
Objective of the E design is to minimize the Euclidean distance~($\|\ve\varphi_1-\ve\varphi_2\|_2$) between the two points (\tikz\draw[line width=0.3 mm, white,fill={rgb,255:red,255; green,100; blue,255}] (0,0) circle (.75ex); in Fig.~\ref{fig:designs}) that belong to the CR and their Euclidean distance is maximal. This can be expressed as
\begin{subequations}\label{eq:oed_nonlin_E}\begin{align}
  &\min_{\stackrel{\ve u_\tau \in [\ve u^L, \ve u^U]}{\forall \tau\in\{\tau_1,\dots,\tau_N\}}}
  \phi_E(\mathcal U)\\
  &\text{s.t. }\phi_E(\mathcal U)=\max_{\ve\varphi_1, \ve\varphi_2} \|\ve\varphi_1-\ve\varphi_2\|_2^2\\
  &\quad\qquad \quad\text{s.t. }\forall j\in\{1,2\}, \ \forall \tau\in\{\tau_1,\dots,\tau_N\}:\\
  &\qquad \qquad\qquad\hat{\ve y}(\ve \varphi_j, \tau) = \ve F_{nl}(\ve \varphi_j, \ve u_{\tau}),\\
  &\quad\qquad \qquad\quad\ve y_m(\tau)=\hat{\ve y}(\hat{\ve p}, \tau),\\
  &\quad\qquad \qquad\quad J_w(\ve\varphi_j) - J_w(\hat{\ve{p}}) \leq \chi^2_{\alpha,n_p}.
\end{align}
\end{subequations}
The E design is also known as a decorrelating design as it aims at finding the experimental conditions such that the CR is as spherical as possible. This criterion is illustrated in Fig.~\ref{fig:designs}. It is noteworthy that the lower-level problem of~\eqref{eq:oed_nonlin_E} scales linearly w.r.t. the number of parameters as it optimizes $2n_p$ variables.

In the linear case the (classical) E design is identified by
\begin{subequations}\label{eq:oed_lin_E}
\begin{align}
  &\min_{\stackrel{\ve u_\tau \in [\ve u^L, \ve u^U]}{\forall \tau\in\{\tau_1,\dots,\tau_N\}}} \tilde\phi_\text{E}(\mathbb{FIM}):=\min_{\stackrel{\ve u_\tau \in [\ve u^L, \ve u^U]}{\forall \tau\in\{\tau_1,\dots,\tau_N\}}} \max_i \lambda_i(\mathbb{FIM}^{-1})\\
  &\qquad\text{s.t. } \hat{\ve{y}}(\ve{p},\tau) = \ve{F}_{l}(\ve{p},\ve{u}_{\tau}),\ \forall \tau \in\{\tau_1,\dots,\tau_N\},
\end{align}
\end{subequations}
where $\lambda_i(\cdot)$ is the $i$-th eigenvalue of a matrix.

\section{Numerical Implementation}\label{sec:implementation}
In this section we discuss the possible ways to solve the proposed optimization problems. {\color{black} We will exploit BARON \citep{baron} in this work in order to guarantee global optimality of the classical OED problems \eqref{eq:oed_lin_A}, \eqref{eq:oed_lin_D}, and~\eqref{eq:oed_lin_E}. Special attention is devoted to the presented bilevel programs as the classical OED problems are single-level optimization problems and can straightforwardly be solved using a nonlinear program solver.} 

We present two {\color{black}simple heuristic} approaches taken from literature that can be used to solve the presented bilevel problems, which can be generalized in the form
\begin{subequations}\label{eq:bilevel_gen}
\begin{align}
\min_{\ve x_1} &\ f(\ve x_1, \ve x_2^*)\\
 \text{s.t. } & \ve x_2^* \in \argmax{\ve x_2} g(\ve x_2)\\
 & \quad \text{s.t. } \ve 0 = \ve h_E(\ve x_1, \ve x_2),\\
 & \quad \ \ \, \quad \ve 0 \geq \ve h_I(\ve x_1, \ve x_2).
\end{align}
\end{subequations}
A special care has to be taken w.r.t. the non-convexity of lower-level problem. Its global optimum has to be identified in order to guarantee feasibility of the upper-level problem~\citep{mit09}.

\subsection{Nested approach}
{\color{black} The following nested approach is inspired by the formulation proposed in~\citep{TANARTKIT1996735,TANARTKIT19971365} for solving dynamic optimization problems and in \citep{RUBEN2013215} for solving a coordination control algorithm using a price-driven coordination technique.}
The nested approach splits the bilevel optimization problem into a lower-level optimization problem
\begin{subequations}\label{eq:bilevel_inn}
\begin{align}
\ve x_2^* \in \ &\argmax{\ve x_2} \ g(\ve x_2)\\
 & \text{s.t. } \ve 0 = \ve h_E(\ve x_1^*, \ve x_2),\\
 & \ \ \, \quad \ve 0 \geq \ve h_I(\ve x_1^*, \ve x_2).
\end{align}
\end{subequations}
that is solved for a given $\ve x_1^*$ using a global solver (e.g. BARON) and a upper-level optimization problem
\begin{subequations}\label{eq:bilevel_gen_b}
\begin{align}
\ve x_1^* \in \ &\argmin{\ve x_1} \ f(\ve x_1, \ve x_2^*)\\
 & \text{s.t. } \ve 0 = \ve h_E(\ve x_1, \ve x_2^*),\\
 & \ \ \, \quad \ve 0 \geq \ve h_I(\ve x_1, \ve x_2^*).
\end{align}
\end{subequations}
that can then be solved for a given $\ve x_2^*$ using a local solver. We use IPOPT~\citep{ipopt} as a local solver in this work.

The individual optimization problems are interconnected by the copy variables $\ve x_1^*$ and $\ve x_2^*$ that are exchanged between the problems. The problems are solved repeatedly and the convergence of the nested approach is claimed once the consecutive values of the copy variables satisfy $\|\ve x^*_{1,k+1} - \ve x^*_{1,k}\|\approx0$ and $\|\ve x^*_{2,k+1} - \ve x^*_{2,k}\|\approx0$, where $k$ is an iteration counter.

If a gradient-based solver is used to determine the local optimum of the problem~\eqref{eq:bilevel_gen_b}, the objective and constraint gradients are to be supplied. An approach from~\cite{Fiacco1990} can be used in this respect with $\ve x_k:=(\ve x_{1,k}^T, \ve x_{2,k}^T)^T$
\begin{align}\label{eq:nested_derivatives}
  \begin{bmatrix} 
  \nabla^2_{\ve x_2,\ve x_2} L|_{\ve x_k} & \nabla_{\ve x_2}^T \ve h|_{\ve x_k}\\
	  - \nabla_{\ve x_2} \ve h|_{\ve x_k} & \ve 0
  \end{bmatrix}  
  \begin{bmatrix} 
    \frac{\text{d}{\ve x_2^*}} {\text{d}{\ve x_1}} \vspace{1mm}\\
    \frac{\text{d}{\ve{\nu}^*}} {\text{d}{\ve x_1}}
  \end{bmatrix} =-
  \begin{bmatrix}
    \nabla^2_{\ve x_2, \ve x_1} L|_{\ve x_k}\\
    \nabla_{\ve x_1} \ve h|_{\ve x_k}
  \end{bmatrix},
\end{align}
where $L$ represents the Lagrangian of the lower-level problem and $\ve\nu:=(\ve\nu_E^T, \ve\nu_{I,i}^T)^T, \forall i\in\mathcal I_A$ is the vector of multipliers corresponding to the equality and active inequality constraints $\ve h(\ve x_1^*, \ve x_2):=(\ve h_E^T(\ve x_1^*, \ve x_2), \ve h_{I,i}^T(\ve x_1^*, \ve x_2))^T, \forall i\in\mathcal I_A$ of the lower-level problem and $\mathcal I_A$ is an index set of the active inequality constraints. {\color{black}It is not guaranteed that the nested approach always converge to a local minimum \citep{TANARTKIT19971365}. Instead it may sometimes converge to a local maximum or a saddle point. The obtained solution can be verified by evaluating the necessary and sufficient conditions for optimality.}

\subsection{KKT-based approach}
Another heuristic approach for solving a bilevel optimization problem is to reformulate the lower-level problem using the KKT conditions~\citep{Dempe2017, StephanDempe, dempe2002}. The reformulated problem reads as
\begin{subequations}\label{eq:bilevel_kkt}
\begin{align}
\min_{\substack{\ve x_1, \ve x_2\\ \ve\nu_E, \ve \nu_I\leq0}} &\ f(\ve x_1, \ve x_2)\\
 \text{s.t. } & \ve 0=\nabla_{\ve x} L(\ve x_1, \ve x_2, \ve\nu_E, \ve\nu_I),\\
 & \ve 0 = \nabla_{\ve\nu} L(\ve x_1, \ve x_2, \ve\nu_E, \ve\nu_I) = \ve h_E(\ve x_1, \ve x_2),\\
 & 0 = \nu_{I,i} h_{I,i}, \ \forall i\in\mathcal I_I,
\end{align}
\end{subequations}
where $\mathcal I_I$ is the index set of the inequality constraints of~\eqref{eq:bilevel_gen}. As discussed above, the lower-level problem has to be solved to global optimality, which is not guaranteed by satisfaction of the KKT conditions. Reaching of global optimum of the lower-level problem has to be assured upon convergence in order to guarantee feasibility of the lower-level problem and thus a local optimum of the bilevel program. The feasibility test can be performed by solving~\eqref{eq:bilevel_inn} globally or by gridding {\color{black}or by set inversion \citep{kieffer2011guaranteed} }techniques with a subsequent comparison of obtained values for variables of lower-level problem.


{\color{black}We note that there are other approaches that can be employed to solve the problem~\eqref{eq:bilevel_gen}. The solution methods proposed by Dutta et al. in \citep{Dutta2006} and by Dempe et al. in \citep{Dempe2016} assume a convex inner level optimization problem. Bard et al., Dempe et al. and Mitsos et al. in \citep{bard1998practical}, \citep{dempe2002}, and \citep{Mitsos2008} proposed solution methods considering a nonconvex inner level optimization problem. It is generally well known that there is a close connection between bilevel problem and semi-infinite programming (SIP) as discussed in \citep{STEIN2002444}. Cutting-plane SIP algorithm is proposed by Blankenship and Falk~\citep{blan1976}, branch and bound algorithm by Bard and Moore~\citep{bard1990} or double penalty function method by Ishizuka and Aiyoshi~\citep{Ishizuka1992}. Recently, Walz et al.~\citep{WALZ201892} presented a global SIP algorithm proposed in~\citep{Mitsos2015} that could be used in the context of optimal experiment design. Reference therein give a complete picture about the use of SIP algorithms for the solution of bilevel programs.}

Also various stochastic methods, such as genetic algorithms, simulated annealing, etc., could be used in principle, where these might be especially interesting for the D design problem~\eqref{eq:oed_nonlin_D} because of its non-smoothness. We only exploit the described nested and KKT-based approac\-hes in this study.

\section{Case Studies}
We apply the presented methodologies for finding OED for two small-scale {\color{black}illustrative} case studies. The employed models are in the form of explicit step responses of linear time-invariant dynamic systems and the optimal experiment design should reveal the best sampling instants. We denote the designs that are based on the exact CRs (problems~\eqref{eq:oed_nonlin_A}, \eqref{eq:oed_nonlin_D}, and \eqref{eq:oed_nonlin_E}) as \emph{exact} designs. The OED problems are solved for $2\sigma$-confidence level ($\alpha = 0.9545$) using the following approaches:
\begin{enumerate}
\item Classical OED problems are solved globally.
\item The designs based on exact CRs and ellipsoidal D design are solved using nested approach where the lower-level problem is solved globally and the upper-level problem is solved using a local solver.
\item In order to study numerical efficacy of different algorithms, the A-optimal design is solved by KKT-based approach globally.
\end{enumerate}

\subsection{Case Study~1}
The mathematical model for biological oxygen demand (BOD)~\citep{bates1988} is considered. The cumulative BOD of microorganisms at incubation time $u_\tau$ is given by
\begin{align}\label{eq:casestudy_BOD}
\hat y(\ve p, \tau) = p_1(1-\exp(-p_2 u_\tau)), \quad u_\tau \in [0,20],
\end{align}
which can also be interpreted as a step response of the first-order linear time-invariant system with static gain $p_1$ and time constant $1/p_2$. At this point it can be observed that $p_1$ enters the output equation linearly while $p_2$ enters nonlinearly.

The true values for the parameters $p_1$ and $p_2$ are, respectively, 2.5 and 0.5. These are also considered as expected least-squares estimates $\hat{\ve p}$ for all the studied OED problems. The measurements $y_m(\tau)$ are assumed to be corrupted by a zero-mean Gaussian noise with the standard deviation 0.1. We will assume here a case where variance or the standard deviation of the measurement noise is unknown. The exact CRs are then defined by~\eqref{eq:NLPE_without_variance}. Additionally, we consider $J(\hat{p})=0$. The tolerance $\epsilon$ for the exact D design is set to $5\times 10^{-3}$.

The optimal sampling times $\mathcal U^*$ with $N=\{4,5\}$ for the classical and exact A-, D- and E-optimal designs and the ellipsoidal D design are reported in Tab.~\ref{tab:meas_times_BOD}. The values of the objective function for exact designs $\phi(\cdot)$ are evaluated at the identified optimal solution $\mathcal U^*$ for each design with $N = \{4,5\}$. For all the designs, the exact OED has a superior performance when compared to the linearized design. The $\mathcal U^*$ as identified by the classical and the exact OED contain multiple common repetitive measurements at $u_\tau^* = 20$. This agreement between the designs can be attributed to the linear entry of $p_1$ into the model $\hat{y}(\tau)$. It can also further be reasoned by an obvious fact; that one can infer the steady-state gain irrespective of the value of time constant closer to the steady state. The classical and the exact A OED had same number of repetitions of $u_\tau^*=20$ using $N=\{4,5\}$, however this is not the case for D-optimal design where only using $N=4$ the number of repetitions match. In case of the E-optimal design, the classical OED identified $\mathcal U^*$ in which $u_\tau^* = 20$ is repeated once more when compared to the solution identified by the exact OED, which again points to the linear decorrelation between $p_1$ and $p_2$ that classical E design tries to achieve. The performance of the proposed D-optimal design based on the inner- and the outer-approximation ellipsoids is better when compared to the classical D design for all values of $N$. It achieves a relatively small loss in performance when compared to the exact D design. This suggests that the proposed ellipsoidal D design is an interesting framework to perform approximate D design as compared to linearization-based alternative.


\begin{table}
\captionsetup{width=\linewidth}
\caption{Optimal designs ($\mathcal U^*$) as identified for classical, ellipsoidal and the exact OED using $N = \{4,5\}$ and the values of objective function of exact designs ($\phi(\cdot)$) evaluated at the identified optimal designs ($\mathcal U^*$). In the case of D design, $\phi(\mathcal U^*):=\hat{\phi}_D(\mathcal{U}^*)$.}
\label{tab:meas_times_BOD}
\centering
\setlength\tabcolsep{3pt} 
\begin{tabular}{c |c c r |c}
 & Design & $N$ & \hspace{1cm} Solution ($\mathcal U^*$) & $\phi(\mathcal U^*)$\\
\toprule[1pt]
\rule{0pt}{3ex}
\multirow{6}{*}{\rotatebox[origin=c]{90}{Classical OED}} 
&A &4 & $\{1.69, 1.69, 20, 20\}$ & 1.610\\ 
&A &5 & $\{1.77, 1.77, 20, 20, 20\}$ & 0.940\\ 
&D &4 & $\{2, 2, 20, 20\}$ & 0.425\\
&D &5 & $\{2, 2, 20, 20, 20\}$ & 0.155\\ 
&E &4 & $\{1.61, 20, 20, 20\}$ & 1.016\\ 
&E &5 & $\{1.75, 20, 20, 20, 20\}$ & 0.365\\ 
\bottomrule[0.5pt]
\rule{0pt}{3ex}
\multirow{2}{*}{\rotatebox[origin=c]{90}{Ellips.}}
&D &4 & $\{1.42, 1.42, 20, 20\}$ & 0.414\\[0.5ex] 
&D &5 & $\{1.69, 1.69, 19.99, 19.99, 20\}$ & 0.154\\ 
\bottomrule[0.5pt]
\rule{0pt}{3ex}
\multirow{6}{*}{\rotatebox[origin=c]{90}{Exact OED}}
&A &4 & $\{1.37, 1.37, 20, 20\}$ & 1.585\\ 
&A &5 & $\{1.60, 1.60, 20, 20, 20\}$ & 0.938\\ 
&D &4 & $\{1.62, 1.62, 20, 20\}$ & 0.409\\ 
&D &5 & $\{1.81, 1.82, 1.83, 19.99, 19.99\}$ & 0.154\\ 
&E &4 & $\{1.04, 1.04, 20, 20\}$ & 0.974\\ 
&E &5 & $\{1.22, 1.23, 20, 20, 20\}$ & 0.322\\ 
\bottomrule[1pt]
\end{tabular}
\end{table}

\begin{figure}[tb!]
\psfrag{p1}[cc][cc][1]{$p_1$}
\psfrag{p2}[cc][cc][1]{$p_2$}
\psfrag{Classical A}[cc][cc][0.8]{\hspace{5mm}Classical A}
\psfrag{Exact A}[cc][cc][0.8]{\hspace{5mm}Exact A}
\centering\includegraphics[width=\linewidth]{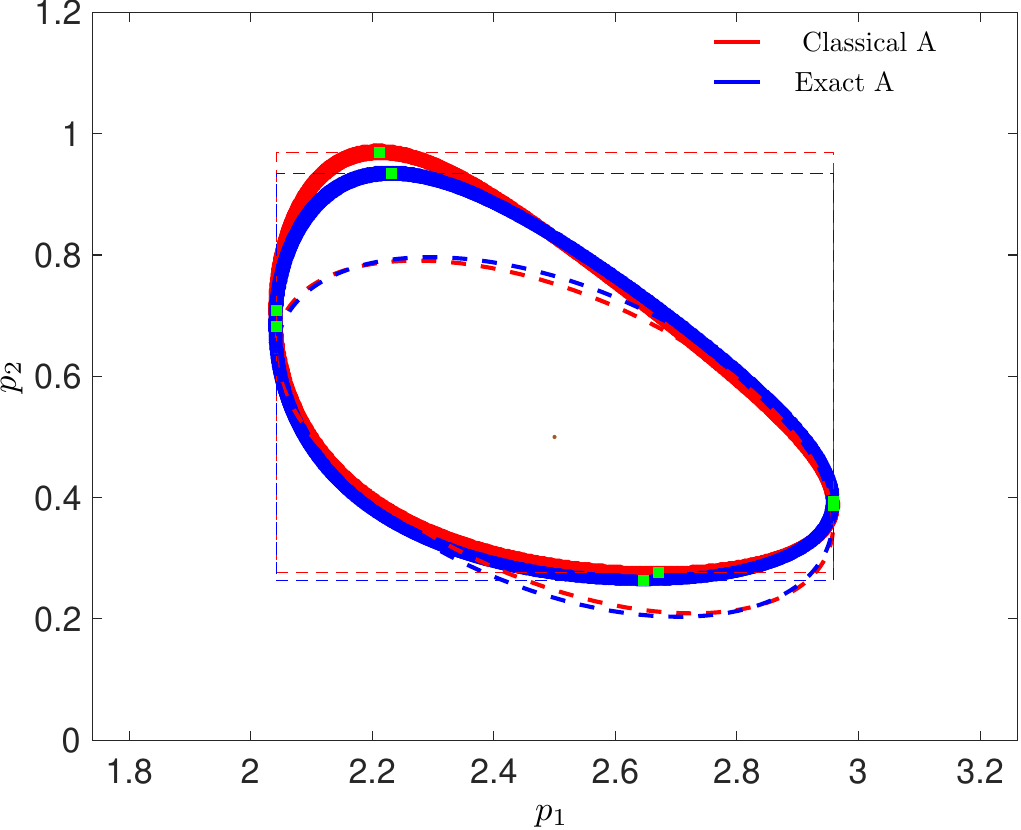}
\caption[]{The exact (solid) and linearized (dashed ellipsoid) CRs using $N = 4$ as obtained by classical and exact A designs. The plot shows the over-approximating orthotopes (dashed) of the exact CRs identified using the anchor points $\ve \pi$ represented by \tikz\draw[line width=0.25 mm,white, fill=green] (0,0) \Square{10pt};.}
\label{fig:A_BOD}
\end{figure}
The exact CRs for the classical (\tikz\draw[line width=0.25 mm, red,fill=red] (0,0) \LongRrectangle{10pt};) and the exact (\tikz\draw[line width=0.25 mm, blue,fill=blue] (0,0) \LongRrectangle{10pt};) A OED for four measurements are compared in Fig.~\ref{fig:A_BOD}. The orthotopes enclosing the exact CRs are plotted using the anchor points $\ve \pi$ represented by 
\tikz\draw[line width=0.25 mm,white, fill=green] (0,0) \Square{10pt};. 
It is clear that the 
orthotope that encloses the exact CR identified by the classical A design (\tikz\draw[dashed,line width=0.25 mm, red,fill=white] (0,0) \Rrectangle{10pt};) has a larger perimeter when one compares it with the orthotope identified by the exact A design. The reason for this can be found when looking at the linearized CRs for both designs (dashed ellipsoids). The linearized ellipsoid clearly does not appro\-xi\-mate the exact CRs well, where, as it can be expected, the approximation is looser w.r.t. $p_2$ that enters nonlinearly in output equation~\eqref{eq:casestudy_BOD}. It is an interesting observation that the presented linearized CRs are very similar to each other, despite the fact that they give significantly different exact CRs.

\begin{figure}[tb!]
\psfrag{p1}[cc][cc][1]{$p_1$}
\psfrag{p2}[cc][cc][1]{$p_2$}
\psfrag{Classical D}[cc][cc][0.8]{\hspace{5mm}Classical D}
\psfrag{Ellipsoidal D}[cc][cc][0.8]{\hspace{6mm}Ellipsoidal D}
\psfrag{Exact D}[cc][cc][0.8]{\hspace{4.5mm}Exact D}
\centering\includegraphics[width=\linewidth]{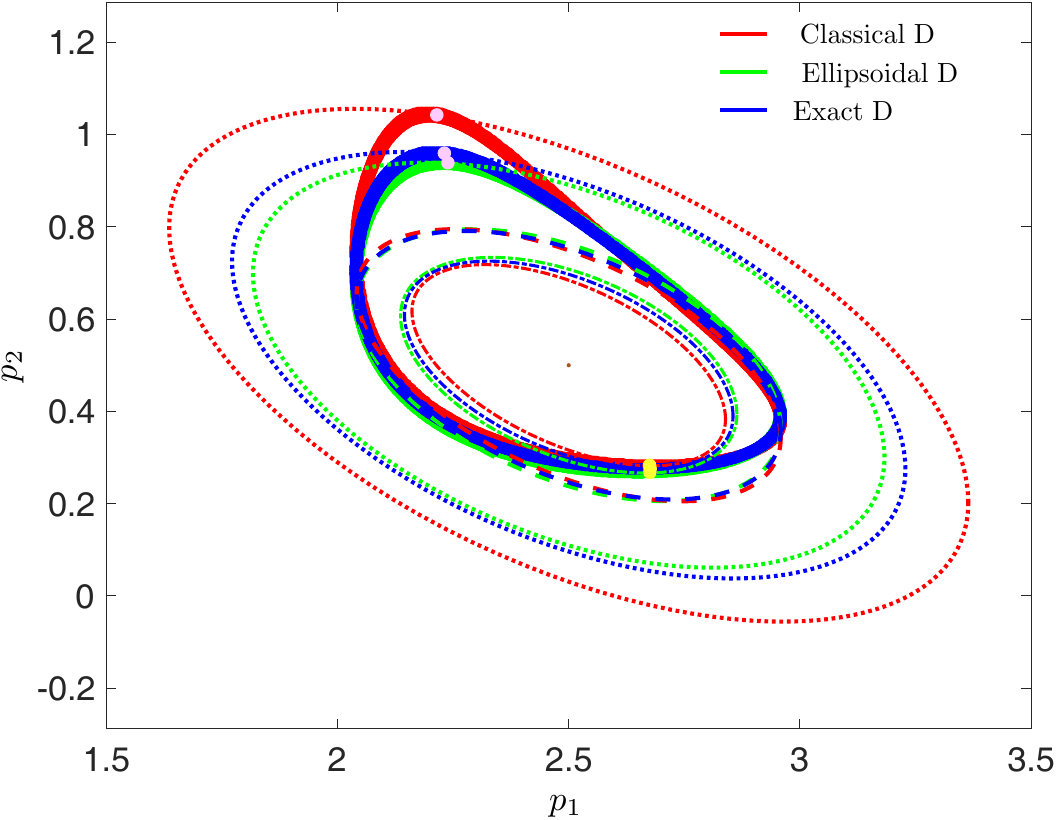}
\caption[]{The exact (solid) and linearized (dashed ellipsoid) CRs using $N = 4$ as obtained by
classical, ellipsoidal and exact D designs. The plot shows the outer-/inner-approximating ellipsoids of the exact CRs (dotted and dash-dotted lines, respectively). \tikz\draw[line width=0.25 mm, white,fill={rgb,255:red,255; green,204; blue,255}] (0,0) circle (.75ex); and \tikz\draw[line width=0.25 mm, white,fill=yellow] (0,0) circle (.75ex); are the intersection points for the outer-/inner-approximating ellipsoids and the exact CRs.}
\label{fig:D_BOD}
\end{figure}
Figure~\ref{fig:D_BOD} shows the exact CRs for the classical (\tikz\draw[line width=0.25 mm, red,fill=red] (0,0) \LongRrectangle{10pt};), ellipsoidal (\tikz\draw[line width=0.25 mm, green,fill=green] (0,0) \LongRrectangle{10pt};) and the exact (\tikz\draw[line width=0.25 mm, blue,fill=blue] (0,0) \LongRrectangle{10pt};) D designs using $N=4$. The linearized CRs (dashed ellipsoids) are again very similar to each other.
Among them the ellipsoid from the classical design (\tikz\draw[dashed,line width=0.4 mm, red,fill=white] (0,0) ellipse (1.75ex and 0.75ex);) has the smallest volume, as might be expected. However, the exact CR for the classical design is the largest one (see Tab.~\ref{tab:meas_times_BOD}), which again comes from disregarding of nonlinearity of the output equation in $p_2$ by the classical design.

We also present the inner- and the outer-appro\-xi\-mation ellipsoids for all the three OED approaches by dash-dotted and dotted ellipsoids, respectively. The corresponding intersection points are marked by 
\tikz\draw[line width=0.25 mm, white,fill=yellow] (0,0) circle (.75ex); and \tikz\draw[line width=0.25 mm, white,fill={rgb,255:red,255; green,204; blue,255}] (0,0) circle (.75ex); 
respectively in Fig.~\ref{fig:D_BOD}. Here we can observe the benefits of weighting introduced in the objective function of the problem~\eqref{eq:oed_nonlin_D_app}. While looking at the sizes of outer-approximating ellipsoids (especially the one constructed for exact D design), it might appear reasonable to minimize the volume of the over-approximating ellipsoid as a good approximation of the exact D design. This would correspond to setting $1/k_{in}\to0$ while solving the problem~\eqref{eq:oed_nonlin_D_app}. We have explored this path in our earlier study~\citep{mupa_ifac17}, but the obtained design results were unsatisfactory, since the over-approximation by an ellipsoid might become very loose. The proposed ellipsoidal D design therefore balances out the concentration on the size of the ellipsoid and the appropriateness of the over-approximation by an ellipsoid.
It is clearly visible that the exact CR for the $\mathcal U^*$ identified by the proposed ellipsoidal OED has much smaller volume than the exact CR identified by the classical approach and, at the same time, it is very close to the optimal exact OED.

\begin{figure}[tb!]
\psfrag{p1}[cc][cc][1]{$p_1$}
\psfrag{p2}[cc][cc][1]{$p_2$}
\psfrag{Classical E}[cc][cc][0.8]{\hspace{5mm}Classical E}
\psfrag{Exact E}[cc][cc][0.8]{\hspace{5mm}Exact E}
\centering\includegraphics[width=\linewidth]{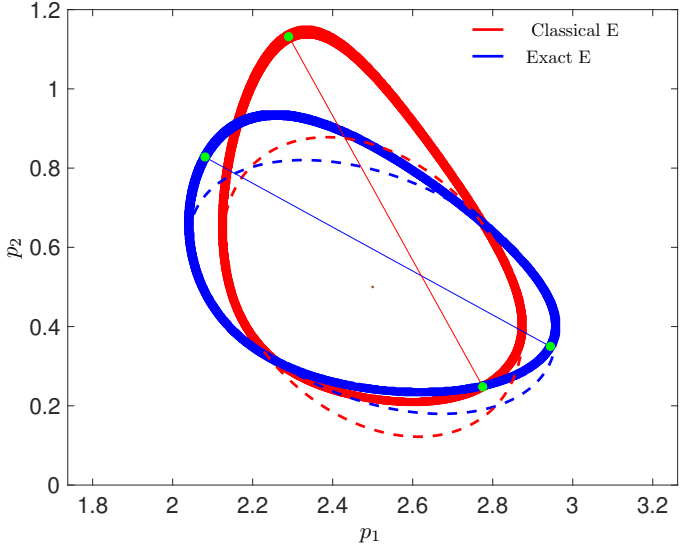}
\caption[]{The exact (solid) and linearized (dashed ellipsoid) CRs using $N = 4$ as obtained for classical and exact E OED. \tikz\draw[line width=0.3 mm, white,fill=green] (0,0) circle (.75ex); mark the points used to calculate the Euclidean distance of the CRs.}
\label{fig:E_BOD}
\end{figure}
The exact CRs for the classical (\tikz\draw[line width=0.25 mm, red,fill=red] (0,0) \LongRrectangle{10pt};) and the exact (\tikz\draw[line width=0.25 mm, blue,fill=blue] (0,0) \LongRrectangle{10pt};) E designs are compared in Fig.~\ref{fig:E_BOD} using $N = 4$ measurements. 
\tikz\draw[line width=0.3 mm, white,fill=green] (0,0) circle (.75ex); in Fig.~\ref{fig:E_BOD} 
mark $\ve\varphi_1$ and $\ve\varphi_2$ (see~\eqref{eq:oed_nonlin_E}) obtained for the classical and the exact E-optimal designs. In this case, unlike for the previous designs, we observe a major discrepancy in the orientation between the linearized ellipsoids obtained for the classical and exact designs, despite the largest semi-axes, which correspond to the largest eigenvalues (see~\eqref{eq:oed_lin_E}) are very similar. Again this is attributed to the nonlinearity in $p_2$. The resulting difference of distances between the most distant points that belong to $P$ is significant, however, among the two designs. Similar behavior can be seen for the case with $N=5$ (see Tab.~\ref{tab:meas_times_BOD}).

\begin{figure}[tb!]
\psfrag{A-design}[cc][cc][1]{A design}
\psfrag{D-design}[cc][cc][1]{D design}
\psfrag{E-design}[cc][cc][1]{E design}
\psfrag{Exact A design Objective}[cc][cc][1]{$\phi_A$}
\psfrag{Exact D design Objective}[cc][cc][1]{$\hat\phi_D$}
\psfrag{Exact E design Objective}[cc][cc][1]{$\phi_E$}
\centering\includegraphics[width=\linewidth]{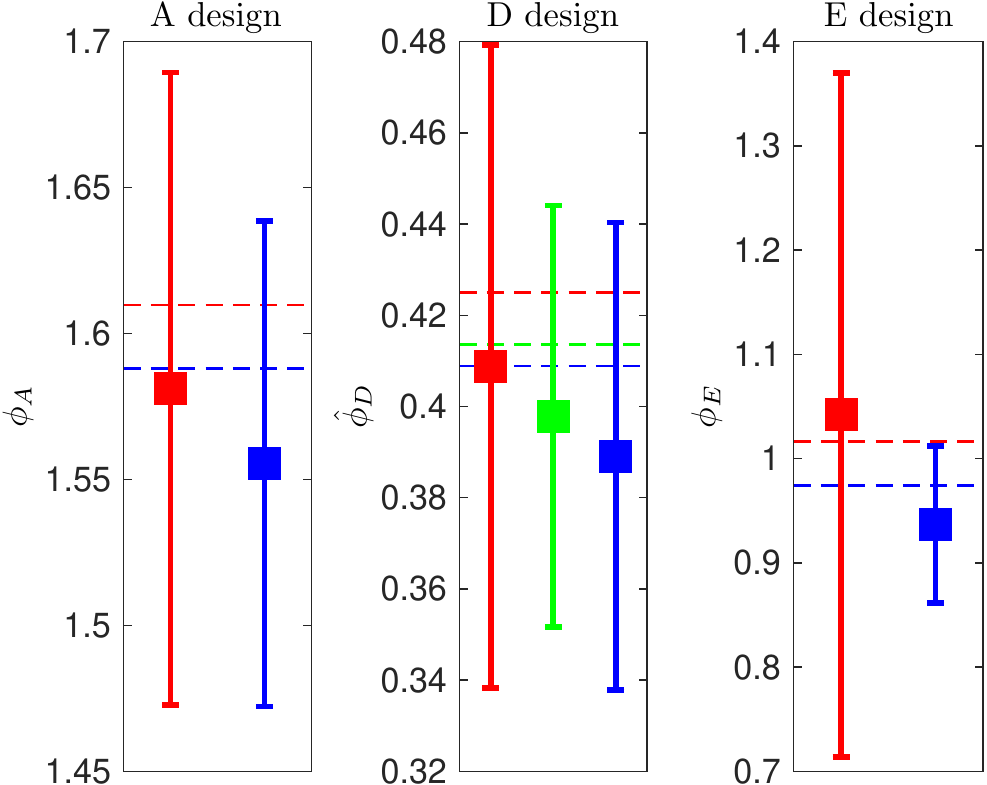}
\caption[]{Mean and variance of $\phi_A$, $\hat\phi_D$ and $\phi_E$ for $1,000$ random experiments with $N=4$ noisy measurements at $\mathcal U^*$ of classical (\tikz\draw[line width=0.25 mm, red,fill=red] (0,0) \LongRrectangle{10pt};), ellipsoidal (\tikz\draw[line width=0.25 mm, green,fill=green] (0,0) \LongRrectangle{10pt};) and the exact (\tikz\draw[line width=0.25 mm, blue,fill=blue] (0,0) \LongRrectangle{10pt};) designs. Dashed line signifies the performance of the nominal design.}
\label{fig:BODr_robustness}
\end{figure}

Next, we study the performance of the obtained designs against a number of simulated experiments. The aim here is to evaluate robustness of the designs against random realization of noise that would be present in the real experiment. We exclude the dependence on the least-squares estimate here, i.e., we will use the nominal values for $\hat{\ve p}$. Such dependence is the subject of study for \emph{robust} OED, which is not in the scope here. We simulated 1,000 experiments with each studied design using the obtained optimal incubation (measurement) times $\mathcal U^*$ and we corrupted the measurements $y_m(\tau) := \hat{y}({\hat{\ve p}},\ve \tau) + e$ with a Gaussian noise $e$ of standard deviation 0.1.

Figure~\ref{fig:BODr_robustness} shows the mean and the variance of the objective value of the exact A, D and E designs. In this plot, we also include the nominal values (when noise-free measurements are gathered) of the different designs using dashed lines. With respect to the mean values of obtained $\phi_A$, $\hat\phi_D$ and $\phi_E$, it is confirmed that exact OED is the best option on average, despite the mean values do not match the expected nominal values of the design criteria. Regarding the obtained variances, it is noteworthy that the classical design exhibits the strongest sensitivity to noise as it can be concluded from the magnitude of the variances and thus appears to be the worst option. This again underpins the importance of consideration of the nonlinearity in the OED and it can be documented by comparing the worst-case value of the exact E design w.r.t. nominal mean obtained for the linearized design (right plot). The last interesting observation is that the variance of $\hat\phi_D$ obtained with ellipsoidal design is slightly smaller than the variance of the exact OED. This might point to the approximation error introduced (tolerance $\epsilon$ in~\eqref{eq:oed_nonlin_D}) in calculating $\hat\phi_D$, which is however not severe as the mean and worst-case values follow the expected trend.

\subsection{Case study~2}
Here we consider a problem where the system output can be modeled as
\begin{equation}\label{eq:cs2_modeleqn}
\hat y(\tau) = b_0\frac{p_1}{p_2^2}\left(\left[p_2\frac{p_1+p_2}{p_1} u_\tau + 1 \right]\exp(-p_2u_\tau) - 1\right),
\end{equation}
which can also be interpreted as a step response of the second-order linear time-invariant system with a double pole at $-p_2$ and a zero at $p_1$. The corresponding transfer function can be given by
\begin{equation}\label{eq:cs2_tfn}
 G(s) = \frac{b_0(s-p_1)}{(s+p_2)^2}.
\end{equation}
Clearly the constant $b_0$ is a parameter that influences the steady-state gain of the system and we will assume it, for simplicity, to be known $b_0=-4$.

Notice that both $p_1$ and $p_2$ enter the output equation nonlinearly, so this problem can be considered as more challenging and even greater discrepancy might be expected between linearized and exact designs. The true values of the parameters $p_1$ and $p_2$ are $0.5$ and $1.0$, respectively, and are equal to $\hat{\ve p}$. The measurements $y_m(\tau)$ are assumed to be corrupted with a zero-mean Gaussian noise with known standard deviation of $0.4$. For a fair comparison of the proposed framework, we assume $J_w(\hat{\ve p})=0$. The tolerance ($\epsilon$) for exact D design is $7.5\times10^{-2}$. 

The classical and the proposed OED frameworks are applied to identify the optimal sampling times $\mathcal U^*$, where $u_\tau^*\in[0,10]$. The previously discussed OED problems are solved for $N = \{2, 3, 4\}$ and $\alpha = 0.9545$~($2\sigma$-CR) with the same numerical techniques as in case study~1. The results are presented in Tab.~\ref{tab:meas_times_second_order_case}. The trends regarding the performance of the different designs are the same as described in case study~1. We observe in this case a more significant benefit of employing an ellipsoidal D design, which, for $N=\{3,4\}$ almost reaches the performance of exact OED and is superior to classical OED.

\begin{table}
\captionsetup{width=\linewidth}
\caption{Optimal designs ($\mathcal U^*$) as identified for classical, ellipsoidal and the exact OED using $N = \{2,3,4\}$  and the values of objective function of exact designs ($\phi(\cdot)$) evaluated at the identified optimal designs ($\mathcal U^*$). In the case of D design, $\phi(\mathcal U^*):=\hat{\phi}_D(\mathcal{U}^*)$.}
\label{tab:meas_times_second_order_case}
\centering
\setlength\tabcolsep{3pt} 
\begin{tabular}{c |c c r |c}
 & Design & $N$ & \hspace{1cm} Solution ($\mathcal U^*$) & $\phi(\mathcal U^*)$\\
\toprule[1pt]
\rule{0pt}{3ex}
\multirow{9}{*}{\rotatebox[origin=c]{90}{Classical OED}} 
&A &2 & $\{1.91,10\}$ & 1.666\\ 
&A &3 & $\{1.86,1.86,10\}$ & 1.151\\ 
&A &4 & $\{1.81,1.81,1.81,10\}$ & 0.974\\ 
&D &2 & $\{2,10\}$ & 0.386\\ 
&D &3 & $\{2,2,10\}$ & 0.231\\ 
&D &4 & $\{2,2,10,10\}$ & 0.148\\ 
&E &2 & $\{1.90,10\}$ & 1.225\\ 
&E &3 & $\{1.82,1.82,10\}$ & 0.520\\ 
&E &4 & $\{1.74,1.74,1.74,10\}$ & 0.341\\ 
\bottomrule[0.5pt]
\rule{0pt}{3ex}
\multirow{3}{*}{\rotatebox[origin=c]{90}{{Ellipsoid\hspace{-1mm}}}}
&D &2 & $\{1.70,10\}$ & 0.363\\ 
&D &3 & $\{1.73,1.73,10\}$ & 0.219\\ 
&D &4 & $\{1.82,1.82,10,10\}$ & 0.144\\ 
\bottomrule[0.5pt]
\rule{0pt}{3ex}
\multirow{9}{*}{\rotatebox[origin=c]{90}{Exact OED}}
&A &2 & $\{1.63,10\}$ & 1.584\\ 
&A &3 & $\{1.67,1.67,10\}$ & 1.132\\ 
&A &4 & $\{1.66,1.66,1.67,10\}$ & 0.966\\ 
&D &2 & $\{1.61,10\}$ & 0.344\\ 
&D &3 & $\{1.65,1.66,10\}$ & 0.218\\ 
&D &4 & $\{1.74,1.77,10,10\}$ & 0.144\\ 
&E &2 & $\{1.62,10\}$ & 1.094\\ 
&E &3 & $\{1.63,1.63,10\}$ & 0.497\\
&E &4 & $\{1.59,1.59,1.59,10\}$ & 0.331\\ 
\bottomrule[1pt]
\end{tabular}
\end{table}

\begin{figure}[tb!]
\psfrag{p1}[cc][cc][1]{$p_1$}
\psfrag{p2}[cc][cc][1]{$p_2$}
\psfrag{Classical A}[cc][cc][0.8]{\hspace{5mm}Classical A}
\psfrag{Exact A}[cc][cc][0.8]{\hspace{5mm}Exact A}
\centering\includegraphics[width=\linewidth]{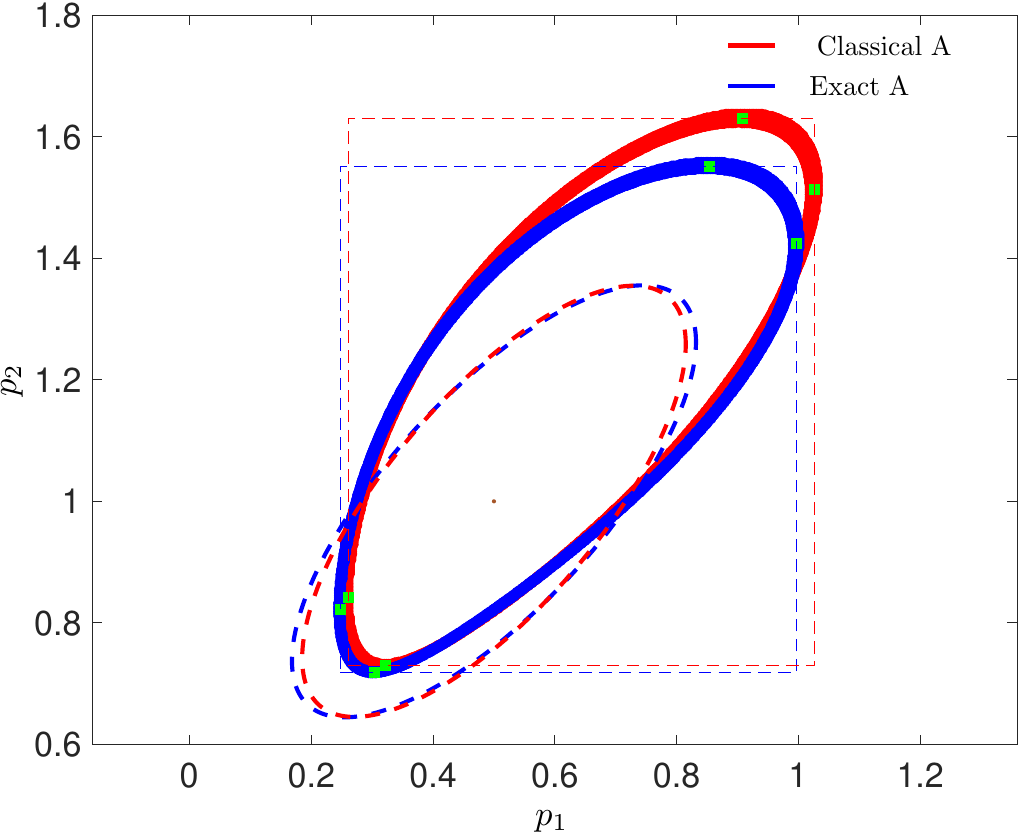}
\caption[]{The exact (solid) and linearized (dashed ellipsoid) CRs using $N = 2$ as obtained for classical and exact A OED. The plot shows the over-approximating orthotopes (dashed) of the exact CRs identified using the anchor points $\ve \pi$ represented by \tikz\draw[line width=0.25 mm,white, fill=green] (0,0) \Square{10pt};.}
\label{fig:A_2ndorder}
\end{figure}
In Fig.~\ref{fig:A_2ndorder}, we compare the exact and linearized CRs obtained for the classical (\tikz\draw[line width=0.25 mm, red,fill=red] (0,0) \LongRrectangle{10pt};) and the exact (\tikz\draw[line width=0.25 mm, blue,fill=blue] (0,0) \LongRrectangle{10pt};) A design using $N = 2$. We can see that the linearized CR captures well the orientation of the exact CR but due to nonlinearity the approximation is relatively poor. Similarly to the first case study we can observe that, despite very similar orientation of the linearized confidence regions, significant benefits of the exact design over the linearized one are obtained. Unlike in case study~1, we obtain reduction in the range on both parametric axes, which is caused by the nonlinearity of the output equation w.r.t. both parameters.

\begin{figure}[tb!]
\psfrag{p1}[cc][cc][1]{$p_1$}
\psfrag{p2}[cc][cc][1]{$p_2$}
\psfrag{Classical D}[cc][cc][0.8]{\hspace{5mm}Classical D}
\psfrag{Ellipsoidal D}[cc][cc][0.8]{\hspace{6mm}Ellipsoidal D}
\psfrag{Exact D}[cc][cc][0.8]{\hspace{4.5mm}Exact D}
\centering\includegraphics[width=\linewidth]{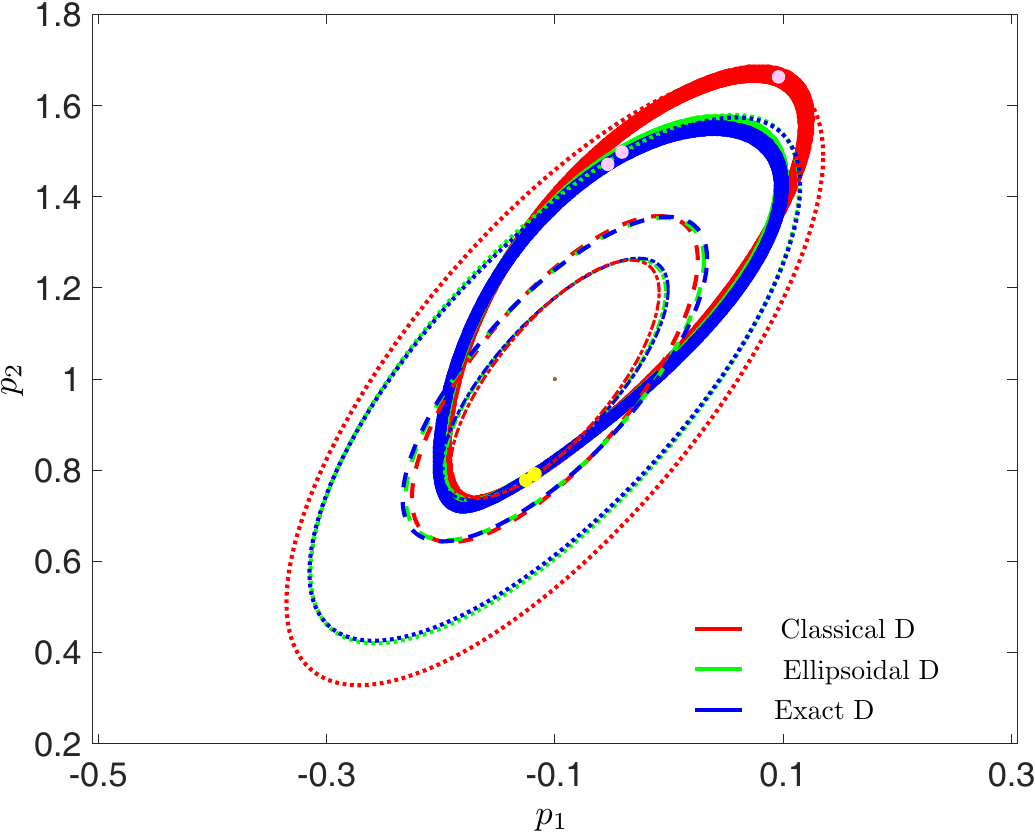}
\caption[]{The exact (solid) and linearized (dashed ellipsoid) CRs using $N = 2$ as obtained for classical, ellipsoidal and exact D OED. The plot shows the outer-/inner-approximating ellipsoids of the exact CRs (dotted and dash-dotted lines, respectively). \tikz\draw[line width=0.25 mm, white,fill={rgb,255:red,255; green,204; blue,255}] (0,0) circle (.75ex); and \tikz\draw[line width=0.25 mm, white,fill=yellow] (0,0) circle (.75ex); are the intersection points for the outer-/inner-approximating ellipsoids and the exact CRs.}
\label{fig:D_2ndorder}
\end{figure}
Figure~\ref{fig:D_2ndorder} shows the resulting CRs for D design criterion. This shows the reason for the good performance of the proposed ellipsoidal technique, which is able to tackle the nonlinearity of the CR far better than the linearization-based design. 

\begin{figure}[tb!]
\psfrag{p1}[cc][cc][1]{$p_1$}
\psfrag{p2}[cc][cc][1]{$p_2$}
\psfrag{Classical E}[cc][cc][0.8]{\hspace{5mm}Classical E}
\psfrag{Exact E}[cc][cc][0.8]{\hspace{5mm}Exact E}
\centering\includegraphics[width=\linewidth]{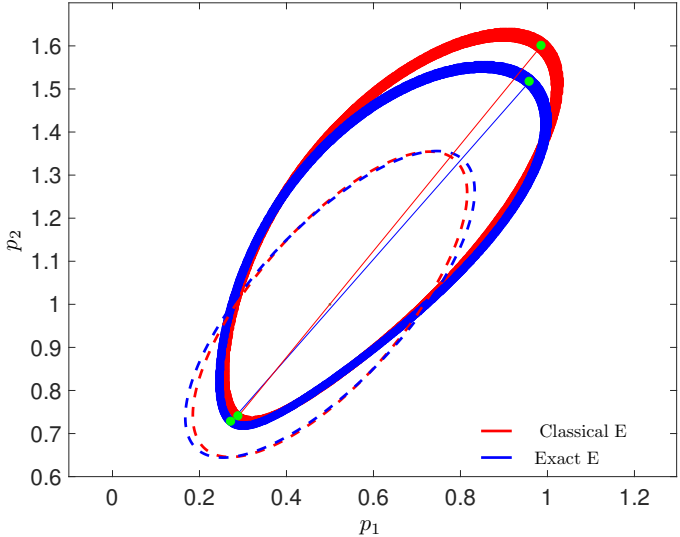}
\caption[]{The exact (solid) and linearized (dashed ellipsoid) CRs using $N = 2$ as obtained for classical and exact E OED. \tikz\draw[line width=0.3 mm, white,fill=green] (0,0) circle (.75ex); mark the points used to calculate the Euclidean distance of the CRs.}
\label{fig:E_2ndorder}
\end{figure}
The E-optimal designs for the classical and the exact OED with $N=2$ are compared in Fig.~\ref{fig:E_2ndorder}. Despite the fact that the linearized CRs show great similarity and they capture the orientation of the exact CR, the exact design tackles the nonlinearity far better and shows clear benefits w.r.t. the linearization-based counterpart.

Robustness of the obtained designs was tested against the random realization of the measurement noise for $N=4$ in the same way as in the previous case study. It is clear again that exact designs outperform the classical OED, which reaches the largest variances and the inferior means. The performance of the ellipsoidal D design is practically the same the performance of exact D design. In comparison with with the case study~1, we observe larger variance of the exact D design, which we attribute to the higher nonlinearity.

\begin{figure}[tb!]
\psfrag{A-design}[cc][cc][1]{A design}
\psfrag{D-design}[cc][cc][1]{D design}
\psfrag{E-design}[cc][cc][1]{E design}
\psfrag{Exact A design Objective}[cc][cc][1]{$\phi_A$}
\psfrag{Exact D design Objective}[cc][cc][1]{$\hat\phi_D$}
\psfrag{Exact E design Objective}[cc][cc][1]{$\phi_E$}
\centering\includegraphics[width=\linewidth]{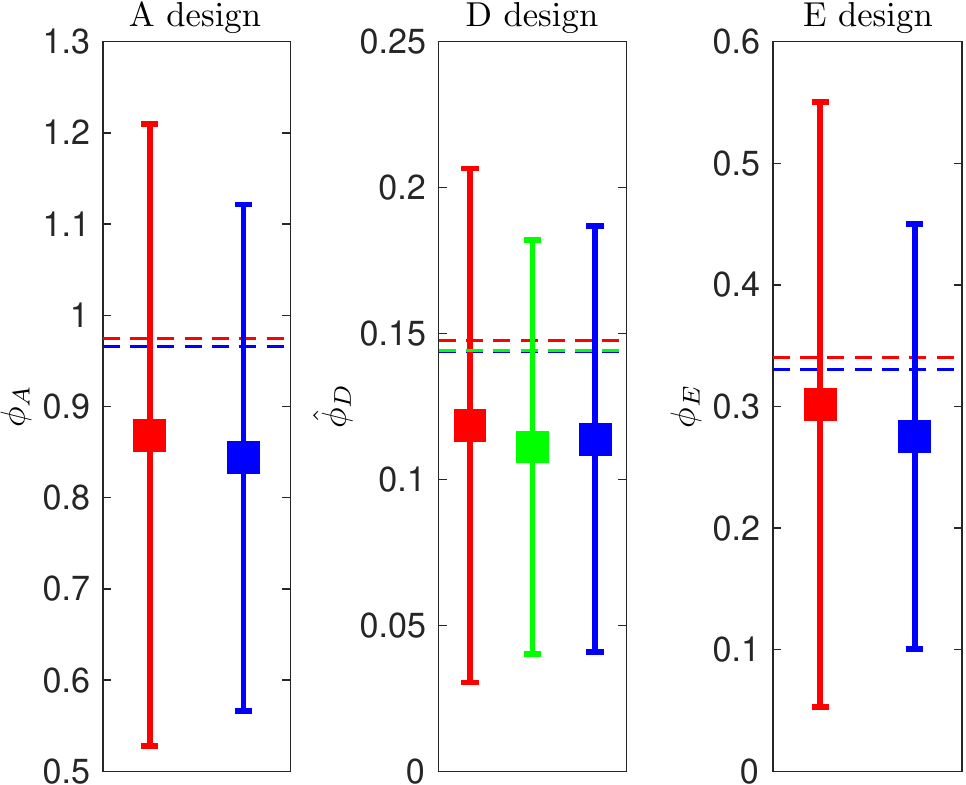}
\caption[]{Mean and variance of $\phi_A$, $\hat\phi_D$ and $\phi_E$ for $1,000$ random experiments using $N=4$ noisy measurements at $\mathcal U^*$ of classical (\tikz\draw[line width=0.25 mm, red,fill=red] (0,0) \LongRrectangle{10pt};), ellipsoidal (\tikz\draw[line width=0.25 mm, green,fill=green] (0,0) \LongRrectangle{10pt};) and the exact (\tikz\draw[line width=0.25 mm, blue,fill=blue] (0,0) \LongRrectangle{10pt};) OED.  Dashed line signifies the performance of the nominal design.}
\label{fig:2ndorder_robustness}
\end{figure}

\subsection{Discussion}
We studied cases where the CR is found for 2$\sigma$ confidence. For a greater confidence, the CR increases in size and nonlinearity affects it stronger. That is why even bigger differences can be expected between classical and exact OED and even greater benefits can be obtained by using exact design (see~\cite{mupa_ifac17}).
 
{\color{black}Regarding the computational efficacy of the different studied problems, it must be clearly pointed out that there exists a high sensitivity of model outputs w.r.t. the model parameters in both case studies due to the presence of highly nonlinear exponential terms. The presence of highly nonlinear terms and the need for inversion of the Fisher information matrix to formulate the objective function in classical OED prohibit BARON from closing the optimality gap and thus it returns locally optimal and sub-optimal solutions unless the optimization problem is properly initialized. However, the classical OED can be solved very efficiently, in few minutes on standard desktop workstation if initialized} {\color{black}properly.} The solution times for exact designs followed the expectations that result from the aforementioned complexity analysis (see Section~\ref{sec:OED}). On average the solution procedure for exact A, D, and E design using nested approach took less than 10\,min, $\approx$6\,h, and 15--30\,min, respectively. This shows that the optimal exact A and E designs can be obtained with practically the same effort as in the case of the classical design for the small-scale problems. The exact A design procedure scales quadratically in $n_p$ so it can get much more time-consuming in higher dimensions. We note that the reduction of CPU time obtained using ellipsoidal D design w.r.t. to exact D design was two-fold ($\approx$3\,h), which is, on one hand, a considerable time saving but, on the other hand, it puts the potential user in question, whether the benefits prevail over the costs. We note, for completeness, that the KKT-based approach applied to problem of exact A design required the solution time of approximately one hour, which makes this approach clearly inferior.

In problems with small number of samples, it might be problematic to identify approximate (experimental) variance or to satisfy the asymptotic properties under which the CRs are defined. In this case, one might think of different approaches to experiment design. One such approach might be to relax the assumption of the presence of a white Gaussian noise in the measurements. This might in turn lead to set-membership estimation approach, also commonly known as guaranteed parameter estimation. A step in the direction of experimental design in set-membership context was taken in~\cite{telen2013guaranteed} and in the recent studies~\citep{WALZ201892, mupa_cace17, mupa_acc16}.

\section{Conclusions}
In this paper, exact and linearization-based methods were presented for the optimal experiment design of a nonlinear parameter estimation problem. The ellipsoidal method is proposed as a computationally less demanding counterpart to the exact D design, which is a computationally demanding problem since it requires a good approximation for the volume of a set. {\color{black}Two simple heuristic} numerical methods are used here to solve the corresponding optimization problems, which are of bilevel nature. The OED framework is tested upon two {\color{black}illustrative} small-scale nonlinear case studies, where the benefits of the exact design are showcased. The proposed ellipsoidal technique is shown to perform very well. Despite this study treated the case when the system model describes a static system in an explicit form, the methodology is straightforwardly applicable to dynamic systems and implicit model equations. An interesting direction for the future work lies, on one hand, in increasing the efficiency of the solution of the bilevel programs and, on the other hand, in the study of robust OED that relaxes the assumption of known (expected) least-squares estimates $\hat{\ve p}$, which might be relevant in practical tasks.

\section*{Acknowledgments}
The research leading to these results has received funding from the European
Commission under grant agreement number 291458 (ERC Advanced Investigator Grant
MOBOCON). RP acknowledges the contribution of the European Commission under the project GuEst (grant agreement No. 790017), of the Slovak Research and Development Agency under the project APVV 15-0007, of the Scientific GrantAgency of the Slovak Republic under the grant 1/0004/17, and of the Research \& Development Operational Programme for the project University Scientific Park STU in Bratislava, ITMS 26240220084, supported by the Research 7 Development Operational Programme funded by the ERDF.

\bibliographystyle{elsarticle-num}

\end{document}